\documentclass[onefignum,onetabnum,final]{siamart171218}%{article}%

\usepackage[utf8]{inputenc} % set input encoding (not needed with XeLaTeX)
\usepackage{graphicx} % support the \includegraphics command and options
\usepackage{multirow}
\usepackage{booktabs} % for much better looking tables
\usepackage{array} % for better arrays (eg matrices) in maths
\usepackage{paralist} % very flexible & customizable lists (eg. enumerate/itemize, etc.)
\usepackage{verbatim} % adds environment for commenting out blocks of text & for better verbatim

\usepackage{epstopdf}

\usepackage{algorithm}
\usepackage{algorithmicx}
\usepackage{algpseudocode}

\usepackage{amsmath,dsfont}
\usepackage{harpoon}
\usepackage{amssymb}
\usepackage{mathrsfs}
\usepackage{graphics}
\usepackage{color}
\usepackage{geometry}
\usepackage{epsfig}
\usepackage{diagbox}
\usepackage{mathtools}
\usepackage{comment}
\usepackage{bbm}
\usepackage{lineno}
\usepackage{fancyhdr}
\usepackage{graphicx}
\usepackage[inline]{changes}
\def\a{\alpha}
\def\b{\beta}

\def\e{\mbox{e}}
\def\dd{\mbox{d}}

\def\g{\gamma}

\def\m{\mathcal}

\def\o{\omega}

\def\t{\textbf}

\def\f{\frac}
\def\p{\partial}

\def\o{\omega}
\def\i{\infty}

\usepackage{wrapfig}
\usepackage{dsfont}
\usepackage{enumerate}
\usepackage{subfigure}
\usepackage{xcolor}
\usepackage{setspace}

\newtheorem{example}{\textup{\textbf{Example}}}

\begin{document}
%\runninglinenumbers
%\allowdisplaybreaks
%\DeclareGraphicsExtensions{.pdf,.jpeg,.png}

\numberwithin{equation}{section} \title{Efficient Scaling and Moving
  Techniques for Spectral Methods in Unbounded Domains}
\author{Mingtao Xia\footnotemark[2] $^,$\footnotemark[3] \and Sihong
  Shao\footnotemark[3] $^,$\footnotemark[1] \and Tom
  Chou\footnotemark[2]}
\renewcommand{\thefootnote}{\fnsymbol{footnote}}

\footnotetext[2]{Department of Mathematics, UCLA,
  Los Angeles, CA, USA 90095-1555}
\footnotetext[3]{LMAM
  and School of Mathematical Sciences, Peking University, Beijing
  100871, CHINA}
\footnotetext[1]{To whom
  correspondence should be addressed. Email:
  \texttt{sihong@math.pku.edu.cn}}
\today
\maketitle

\begin{abstract}
When using Laguerre and Hermite spectral methods to numerically solve
PDEs in unbounded domains, the number of collocation points assigned
inside the region of interest is often insufficient, particularly when
the region is expanded or translated to safely capture the unknown
solution.  Simply increasing the number of collocation points cannot
ensure a fast convergence to spectral accuracy. In this paper, we
propose a scaling technique and a moving technique to adaptively
cluster enough collocation points in a region of interest in order to
achieve a fast spectral convergence.  Our scaling algorithm employs an
indicator in the frequency domain that is used to determine when
scaling is needed and informs the tuning of a scaling factor to
redistribute collocation points to adapt to the diffusive behavior of
the solution.  Our moving technique adopts an exterior-error indicator
and moves the collocation points to capture the translation.  Both
frequency and exterior-error indicators are defined using only the
numerical solutions. We apply our methods to a number of different
models, including diffusive and moving Fermi-Dirac distributions and
nonlinear Dirac solitary waves, and demonstrate recovery of spectral
convergence for time-dependent simulations.    
Performance comparison in solving a linear parabolic problem
shows that our frequency scaling algorithm outperforms the existing scaling approaches.
We
also show our frequency scaling technique is able to track the blowup
of average cell sizes in a model for cell proliferation.

%Comparing the performance
%of our frequency scaling algorithm while solving a linear parabolic
%problem shows that it outperforms the existing scaling approaches.

%
%including interpolation of diffusing and
%moving Fermi-Dirac distributions and nonlinear Dirac solitary
%waves, and a comparison study on the parabolic model, demonstrate
%their ability in recovering the spectral convergence for
%time-dependent simulations.  
%

\vspace*{4mm}

\end{abstract}

\section{Introduction}

Many scientific models described by PDEs with blowup solutions are set
in unbounded domains.  For example, in many models of cellular
proliferation, a ``blowup'' in which the average size of a population
of cells becomes uncontrolled and diverges over many generations of
growth is possible \cite{KESSLER2017}. The conditions under which
blowup occurs is difficult to determine analytically
\cite{BERNARD2016} but has been explored numerically \cite{Xia2020}.
However, numerically tracking ``blowup'' behavior over long times is
extremely difficult, as it requires solving the problem in a truly
unbounded domain to capture the diverging mean size.  There are many
other problems where it is desirable to find a numerical solution in
an unbounded domain, including the stability of solitary waves arising
from the nonlinear Dirac equation \cite{Shao2014Stability,Cuevas2016} and
diffusion in a parabolic system \cite{MA2005}.

Considerable progress has recently been made in spectral methods for
solving PDEs in unbounded domains \cite{shen2009some}.  Among the
existing spectral methods, the direct approach that is typically used
is based on orthogonal basis functions defined on infinite intervals,
\textit{e.g.}, the Hermite and Laguerre spectral methods
\cite{COULAUD1990,GUO2006,Tao2018Hermite}.  It has been demonstrated
that the performance of these spectral methods can be greatly improved
when a proper coordinate scaling is used
\cite{tang1993hermite,shen2009some}.  However, it is not clear how to
systematically perform the scaling, especially when transient behavior
arises.  A Hermite spectral method with time dependent scaling has
been proposed for parabolic problems by introducing a time dependent
scaling factor $\b(t)$ to meet the coercive condition \cite{MA2005}.
Nonetheless, the form of $\b(t)$ and related parameters are chosen
based on specified knowledge of parabolic models and thus cannot be
easily generalized to other problems.

Motivated by the success of adaptive methods in bounded domains
\cite{ren2000iterative,Tang2003ADAPTIVE,li2001moving}, we propose two
indicators to adaptively allocate a sufficient number of collocation
points to represent the unknown solution in the region of
interest. The first indicator, designed for matching the diffusion of
unknown solutions, extracts the frequency-space information of
intermediate numerical solutions and isolates its high frequency
components. This frequency indicator not only provides a lower bound
for the interpolation error, but also measures the decay of the
derivatives of the reference solution as $|x|\to +\infty$. By tuning a
scaling factor in our proposed scaling technique, the frequency
indicator can be maintained at a low level. However, the translation
of unknown solutions may also amplify the frequency indicator and thus
may result in larger errors for excessive scaling. To accommodate this
scenario, a second, exterior-error indicator is used to calculate an
upper bound for the error in the exterior domain, allowing one to
capture translation via moving collocation points.  Accordingly, for
problems that may involve both translation and diffusion in unbounded
domains, the above two indicators are combined in a ``first moving
then scaling'' approach. Numerical experiments demonstrate their
ability to recover a faster spectral convergence for time-dependent
solutions.

%However, when solving the unbounded domain problem with spectral
%methods directly without compromising to a bounded domain, proper
%scaling is required as stated in \cite{tang1993hermite} to
%redistribute the collocation points and match diffusive behavior
%especially in a time dependent problem. Furthermore, moving the basis
%functions to match translative behavior and introduce new collocation
%points are also necessary. Therefore, in this paper we aim at
%devising techniques that could be applied to adaptively adjust
%collocation points over time and ensure that there are always enough
%collocation points in regions of interests and numerical accuracy can
%be maintained.

The remainder of this paper is organized as follows.
Section~\ref{sec:scaling} introduces the frequency indicator,
connects it to the approximation
error,
and proposes the frequency-dependent scaling
technique for diffusion. 
Section~\ref{sec:moving} proposes the exterior-error-dependent
moving technique for translating problems. We then combine, in
Section~\ref{sec:scalingmoving}, the above two approaches to solve
time-dependent problems involving both diffusion and
translation. Section~\ref{sec:com} compares the frequency-dependent
scaling with a time-dependent scaling proposed in \cite{MA2005} for
solving parabolic systems. In Section~\ref{sec:bio}, we apply the
frequency-dependent scaling method to a PDE model describing
structured cell populations to track blowup behavior. Finally, we
summarize our approaches and make concluding remarks in
Section~\ref{sec:con}.

%%%%%%%%%%%%%%%%%%%%%%%%%%%%%%%%%%%%%%%%%%%%%%%%%%%%%%%%%%%%%%%%%%%%%%%%%%%%
%%%%%%%%%%%%%%%%%%%%%%%%%%%%%%%%%%%%%%%%%%%%%%%%%%%%%%%%%%%%%%%%%%%%%%%%%%%%
\section{Frequency-dependent scaling}
\label{sec:scaling}

In this section, we formulate a scaling technique by first extracting
frequency domain information on the evolution of numerical
solutions, the pseudo-code of which is presented in
Alg.~\ref{algscaling}. Following Guo \textit{et al.} \cite{GUO2006},
the discussion utilizes the generalized Laguerre polynomials which are
mutually orthogonal on the half-line $\Lambda\coloneqq(0, +\i)$ with weight
function

\begin{equation}
\o_{\a, \b}(x)=x^{\a}\e^{-\b{x}}, \quad \a >-1,\, \b > 0.
\end{equation}
The generalized Laguerre polynomials of degree $\ell$ are denoted by
$\m{L}_{\ell}^{(\a, \b)}(x)$ and reduce to the usual Laguerre
polynomials when $\beta=1$.  In this work, we regard $\beta$ to be the
{\em scaling factor}, and seek a time-dependent spectral approximation
of $u(x,t)$ on $\Lambda$. Henceforth, 
for notational simplicity, the $t$-dependence will usually be omitted.

%The $t$-dependence is will usually be
%neglected for notational simplicity.

For any $u\in L^2_{\o_{\a, \b}}(\Lambda)$, the spectral
approximation using the interpolation operator $\mathcal{I}_{N,
  {\a,\b}}$ is

\begin{equation}
\label{SpecD}
u(x) \approx U_N^{(\a, \b)}(x) = \mathcal{I}_{N, {\a, \b}}u
=\sum_{\ell=0}^{N}u_{\ell}^{(\a, \b)}\m{L}_{\ell}^{(\a, \b)}(x),
\end{equation}
where the coefficients $u_{\ell}^{(\a, \b)}$ can be computed by using
\textit{e.g.}, the Laguerre-Gauss collocation points $x_j^{(\a, \b)}$,

\begin{equation}\label{eq:coef}
  {u}^{(\a, \b)}_{\ell} = \f{1}{\g^{(\a, \b)}_{\ell}}
  \sum_{j=0}^N\m{L}_{\ell}^{(\a, \b)}(x_j^{(\a, \b)})u(x_j^{(\a, \b)})w_j^{(\a, \b)},
\quad \ell=0,1,\ldots,N,
\end{equation}
where $N$ is the expansion order (\textit{i.e.}, $N+1$ collocation
points or $N+1$ basis functions), $\g_{\ell}^{(\a,
  \b)}=(\m{L}_{\ell}^{(\a, \b)},\m{L}_{\ell}^{(\a, \b)})_{{\o_{\a,
      \b}}}$ is the $L^2_{\o_{\a, \b}}$ inner product, $w_j^{(\a,
  \b)}$ denotes the corresponding weight for collocation point
$x_{j}^{(\alpha, \beta)}$, and

\begin{equation}
  u(x_j^{(\a, \b)}) = U_N^{(\a, \b)}(x_j^{(\a, \b)})
  =\mathcal{I}_{N, {\a, \b}}u(x_j^{(\a, \b)}),
\quad j=0,1,\ldots, N.
\label{NUMSOL}
\end{equation}

Let $A^r_{\a, \b}(\Lambda)$ be the nonuniformly weighted Sobolev space.
For any integer $r \ge 0$, its seminorm and norm are defined by

\begin{equation}
|u|_{A^r_{\a,\b}} = \|{\partial_x^r u}\|_{\o_{\a+r,\b}}, \quad
\|u\|_{A^r_{\a,\b}} = \left(\sum_{k=0}^r |u|_{A^k_{\a,\b}}^2 \right)^{1/2}.
\end{equation}

For any $u\in A^r_{\a-1, \b}(\Lambda)\cap A^r_{\a, \b}(\Lambda)$ with
integer $r\ge 1$, there is a well-known interpolation error estimate
when using Laguerre-Gauss collocation points \cite{GUO2006}:

 \begin{equation}
  \|\m{I}_{N, {\a, \b}}u - u\|_{\o_{\a, \b}}\leq
  c(\b N)^{\f{1-r}{2}}(\b^{-1}|u|_{A_{\a-1, \b}^r}+(1+\b^{-\f{1}{2}})
  (\mbox{ln}N)^{\f{1}{2}}|u|_{A_{\a, \b}^r}).
\label{ESTIMATE}
\end{equation}
Here, $c$ denotes a generic positive constant which does not depend on
$\a$, $\b$, $N$, or any function. This error estimate is a crucial
element in the formal development and successful implementation of the
proposed scaling and moving techniques.

When the scaling factor is updated from $\beta$ to $\tilde{\beta}$,
the collocation points, weights and $L^2_{\o_{\a, \b}}$ norms are
updated according to

\begin{equation}
x_j^{(\a, \tilde{\b})} = \f{\b}{\tilde{\b}} x_j^{(\a, \b)}, \quad
w_j^{(\a, \tilde{\b})}  = \f{\b^{\a+1}}{\tilde{\b}^{\a+1}} w_j^{(\a, \b)}, \quad
\g^{(\a, \tilde{\b})}_{\ell}  = \f{\b^{\a+1}}{\tilde{\b}^{\a+1}} \g^{(\a, \b)}_{\ell}.
\end{equation}
The expansion coefficients $u^{(\a, \tilde{\b})}_{\ell}$ can then be
estimated through Eq.~\eqref{eq:coef} where we may use the
approximation \eqref{SpecD}: $u(x_j^{(\a, \tilde{\b})})\approx
U_N^{(\a, \b)}(x_j^{(\a, \tilde{\b})})$. This procedure constitutes
the $\Call{scale}{}$ subroutine in Lines~\ref{alg:s0} and \ref{alg:s1}
of Alg.~\ref{algscaling}.

To implement the scaling technique, one needs to determine when to
apply it and how to choose a new scaling factor $\tilde{\beta}$ such
that spectral accuracy can be kept for a prescribed expansion of order
$N$.  To this end, we propose a {\em frequency indicator} acting on
the numerical solution $U_N^{(\a, \b)}$:

\begin{equation}
\m{F}(U^{(\a, \b)}_N) =
\left({\f{\sum\limits_{\ell=N-M+1}^{N}\g_{\ell}^{(\a, \b)} (u^{(\a,
      \b)}_{\ell})^2}{\sum\limits_{\ell=0}^{N}\g_{\ell}^{(\a,
      \b)}(u^{(\a, \b)}_{\ell})^2}}\right)^{\f{1}{2}},
\label{FREQSCAL}
\end{equation}
which measures the contribution of the $M$ highest-frequency
components to the $L^2_{\o_{\a, \b}}$-norm of $U_N^{(\a, \b)}$.  The
subroutine $\Call{frequency\_indicator}{}$ in Lines~\ref{alg:i00},
\ref{alg:i0}, \ref{alg:i1}, and \ref{alg:i2} of Alg.~\ref{algscaling}
calculates this contribution in which we choose $M = [\f{N}{3}]$ in
view of the often-used $\frac{2}{3}$-rule \cite{Hou2007Computing,
  Orszag1971On}.

If the frequency indicator $\m{F}(U_N^{(\a, \b)})$ increases over
time, the contribution of high frequency components to the numerical
solution increases, indicating that the numerical solution is decaying
more slowly in $x$ and that we need to adjust the scaling factor to
enlarge the computational domain $[x_0^{(\a, \b)}, x_N^{(\a, \b)}]$
demarcated by the smallest and largest collocation point positions.
In Line \ref{alg:cond} of Alg.~\ref{algscaling}, $\nu f_0$ is the
threshold at some time $t$. If the value of the frequency indicator of
the current numerical solution $f>\nu f_0$, then we consider scaling.
The parameter $\nu$ is usually chosen to be slightly larger than $1$
to prevent the frequency indicator becoming too large without invoking
scaling.

However, the \t{if} condition is only a necessary condition. Only
after we enter the \t{while} loop in Line \ref{alg:scalewhile} will we
perform scaling, which aims to ensure that the frequency indicator
$\m{F}(U^{(\a, \b)}_N)$ will not increase after scaling. Actually,
this \t{while} loop tries to minimize $\m{F}(U^{(\a, \b)}_N)$ by
geometrically shrinking the scaling factor $\beta$ ($q$ in Line
\ref{q} is the common ratio) to ensure sufficient scaling since
$\m{F}(U^{(\a, \b)}_N)$ is a lower bound for the numerical error, as
shown in Eq.~\eqref{lower}.  A more continuous adjustment is preferred
by setting $q$ to be slightly less than $1$, which may also prevent
over-shrinking of the scaling factor within one single time step.
Henceforth, we will choose $q=0.95$ and $\nu=1/q$.  Moreover, at the
initial time $t=0$, we also ensure the frequency indicator is small
enough by choosing a suitable initial scaling factor.

\begin{algorithm}[t]
\caption{\small Pseudo-code of spectral methods with
  frequency-dependent scaling.}
\begin{algorithmic}[1]
\State Initialize $N$, $\nu>1$, $q<1$, $\Delta t$, $T$, $\alpha$, $\beta$, $U_{N}^{(\alpha,\beta)}(0)$, $\underline{\b}$
\State $t \gets 0$
\State $f_0 \gets \Call{frequency\_indicator}{U_{N}^{(\alpha,\beta)}(t)}$ \label{alg:i00}
\While{$t<T$}
\State $U_{N}^{(\alpha,\beta)}(t+\Delta t)\gets \Call{evolve}{U_{N}^{(\alpha,\beta)}(t),\Delta t}$
\State $f\gets \Call{frequency\_indicator}{U_{N}^{(\alpha,\beta)}(t+\Delta t)}$ \label{alg:i0}
\If{$f > \nu f_0$} \label{alg:cond}
\State $\tilde{\b}  \gets q \beta$
\State $U_N^{(\a,\tilde{\b})}  \gets  \Call{scale}{U_{N}^{(\alpha,\beta)}(t+\Delta t),\tilde{\b}}$ \label{alg:s0}
\State $\tilde{f} \gets \Call{frequency\_indicator}{U_N^{(\a,\tilde{\b})}}$ \label{alg:i1}
\While{$\tilde{f}\leq f$ and $\tilde{\b}\geq \underline{\b}$}\label{alg:scalewhile}
\State $\b \gets \tilde{\b}$
\State $U_{N}^{(\alpha,\beta)}(t+\Delta t)\gets U_N^{(\a,\tilde{\b})}$
\State $f_0 \gets \tilde{f}$
\State $f \gets \tilde{f}$
\State $\tilde{\b}  \gets q \beta$ \label{q}
\State $U_N^{(\a,\tilde{\b})}  \gets  \Call{scale}{U_{N}^{(\alpha,\beta)}(t+\Delta t),\tilde{\b}}$  \label{alg:s1}
\State $\tilde{f} \gets \Call{frequency\_indicator}{U_N^{(\a,\tilde{\b})}}$ \label{alg:i2}
\EndWhile
\EndIf
\State $t\gets t+\Delta t$
\EndWhile
\end{algorithmic}\label{algscaling}
\end{algorithm}

In this work, the generalized Laguerre polynomials with $\alpha=0$ are
used and the relative $L^2_{\o_{\a, \b}}$-error

\begin{equation}
\text{Error} = \f{\|U^{(\a, \b)}_N-u\|_{\o_{\a, \b}}}{\|u\|_{\o_{\a, \b}}},
\label{RelaE}
\end{equation}
is used to measure the quality of the spectral approximation $U^{(\a, \b)}_N(x)$
to the reference solution $u(x)$.  We always use the most updated
scaling factor to calculate the above error.

 \begin{figure}[htb]
\begin{center}
    \includegraphics[width=6.2in]{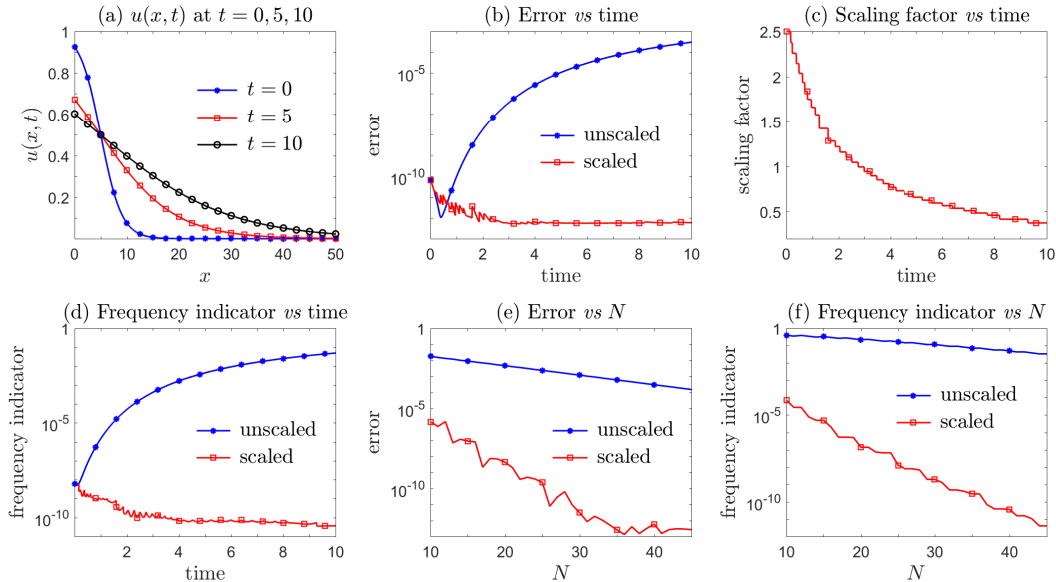}
\end{center}
\vspace{-3mm}
     \caption{\small Numerical approximation to the diffusive
       Fermi-Dirac distribution $u(x,t)$ given by Eq.~\eqref{eq:fd}.
       The scaling algorithm~\ref{algscaling} produces much more
       accurate solutions and recovers a faster spectral convergence with
       respect to the expansion order $N$. As we
       expected, the frequency indicator defined in
       Eq.~\eqref{FREQSCAL} shows a similar behavior to the error
       defined in Eq.~\eqref{RelaE} against either time or $N$. The
       data in last two plots are measured at $t=10$. }
     \label{Fig1}
\end{figure}

\begin{example}\label{ex:fd}
\rm
We use the spreading Fermi-Dirac distribution

\begin{equation}\label{eq:fd}
u(x, t) = \f{1}{1+\e^{\f{x-5}{2+t}}},
\end{equation}
to test the performance of the scaling algorithm~\ref{algscaling}.  It can be readily verified that the reference
solution $u(x, t)$ expands over time as shown in
Fig.~\ref{Fig1}(a). The proposed frequency-dependent scaling with
$N=40$ effectively maintains the relative error under $10^{-10}$ up
until time $t=10$ whereas the error for the corresponding unscaled
solution rapidly grows to over $10^{-4}$ (see Fig.~\ref{Fig1}(b)). We
also plot, as $u(x,t)$ evolves, the history of the scaling factor
$\beta$ and frequency indicator $\m{F}(U^{(\a, \b)}_N)$ in
Figs.~\ref{Fig1}(c) and \ref{Fig1}(d), respectively.  It is clear that
the frequency indicator increases for the unscaled solution as time
evolves and that time-dependent scaling is required to preserve the
accuracy.  The proposed frequency-dependent scaling technique detects
the error and shrinks the scaling factor in order to enlarge the
computational domain in accordance with the expansion of the reference
solution. The spectral convergence as a function of the expansion
order $N$ can be also recovered by Alg.~\ref{algscaling}. The errors at the final time, for the scaled
and unscaled approach, are displayed in Fig.~\ref{Fig1}(e).  The final
scaling factors at $t=10$ are $0.3213, 0.3560, 0.3747, 0.3945, 0.3945$
for $N=25,30,35,40,45$, respectively, having all decreased from the
common initial scaling factor of 2.5. Figs.~\ref{Fig1}(e, f) show very
similar and expected behavior of the frequency indicator and error as
a function of $N$. Since the error and the frequency indicators behave
similarly across time (see Figs.~\ref{Fig1}(b, d)), we also expect them
to behave similarly with $N$.  These similarities suggest a possible
connection between the error and the frequency indicator.
  
\end{example}

The success of the scaling algorithm~\ref{algscaling} is rooted
in the connection between the frequency indicator \eqref{FREQSCAL} and
the evolution of the information embedded in the numerical solutions.
There are two reasons to use a frequency indicator.  First, starting
from Eq.~\eqref{FREQSCAL} with $M = [\f{n}{3}]$ and a sufficiently
large expansion order $N$, we have
\begin{align}
  \frac{1}{2} \m{F}(U^{(\a, \b)}_N) &\approx \frac{1}{2}
  \frac{\|\m{I}_{N, {\a,\b}}u-\m{I}_{N-M, {\a,\b}}u\|_{\o_{\a,\b}}}
       {\|\m{I}_{N, {\a,\b}}u\|_{\o_{\a,\b}}} \nonumber \\
\: & \le \frac{1}{2} \f{\|u - \m{I}_{N, {\a, \b}}u\|_{\o_{\a, \b}}+ \|u -
  \m{I}_{N-M, {\a, \b}}u\|_{\o_{\a, \b}}}{\|\m{I}_{N, {\a, \b}}u\|_{\o_{\a, \b}}}
\label{lower}\\
\: & \le \f{\|u -\m{I}_{N-M, {\a, \b}}u\|_{\o_{\a, \b}}}
   {\|\m{I}_{N, {\a, \b}}u\|_{\o_{\a, \b}}}, \nonumber
\end{align}
which provides an estimate to the lower bound of $\|u - \m{I}_{N-M,
  {\a, \b}}u\|_{\o_{\a,\b}}$. Minimizing $\m{F}(U^{(\a, \b)}_N)$ in
Alg.~\ref{algscaling} may reduce the lower bound of the interpolation
error. Moreover, a straightforward application of the interpolation
error estimator \eqref{ESTIMATE} to the two terms in the numerator of
Eq.~\eqref{lower} yields

\begin{equation}
  \label{ieq:two}
  \left(\sum\limits_{\ell=N-M+1}^{N}
  \g_{\ell}^{(\a, \b)}(u^{(\a, \b)}_{\ell})^2\right)^{1/2}
  \leq c_F(\b N)^{\f{1-r}{2}}\left(\b^{-1}|u|_{A_{\a-1, \b}^r}
  +(1+\b^{-\f{1}{2}})(\mbox{ln}N)^{\f{1}{2}}|u|_{A_{\a, \b}^r}\right),
\end{equation}
where the constant $c_F \equiv (1 + 2^{\f{r-1}{2}})c$.  Thus, we
find

\begin{equation}\label{ieq:estF}
  \m{F}(U^{(\a, \b)}_N)\le c_F(\b N)^{\f{1-r}{2}}
  \left(\b^{-1}\f{|u|_{A_{\a-1, \b}^r}}{ \|U^{(\a, \b)}_N\|_{\o_{\a,\b}}}
  +(1+\b^{-\f{1}{2}})(\mbox{ln}N)^{\f{1}{2}}\f{|u|_{A_{\a, \b}^r}}{\|U^{(\a, \b)}_N\|_{\o_{\a,\b}}}\right),
\end{equation}
implying that $\forall\, \varepsilon\in(0,1)$, we may choose a
sufficiently large $N$ such that $\m{F}(U^{(\a, \b)}_N)<\varepsilon$.

Secondly, the frequency indicator $\m{F}(U^{(\a, \b)}_N)$ can be used
to measure the decay of the reference solution's derivatives as $x$
tends to infinity. According to inequality \eqref{ieq:estF}, if
${|u|_{A_{\a-1, \b}^r}}/{\|U^{(\a, \b)}_N\|_{\o_{\a,\b}}}$ is fixed, a
larger $\m{F}(U^{(\a, \b)}_N)$ implies a larger ${|u|_{A_{\a,
      \b}^r}}/{\|U^{(\a, \b)}_N\|_{\o_{\a,\b}}}$.  In particular,
given $s\in\Lambda$ (\textit{e.g.}, $s =\sqrt{2x_N^{(\a, \b)}}$), if

\begin{equation}
\m{F}(U^{(\a, \b)}_N)> c_F (\b N)^{\f{1-r}{2}}\f{|u|_{A_{\a-1, \b}^r}}{\|U^{(\a, \b)}_N\|_{\o_{\a,\b}}}(\b^{-1} + s (1+\b^{-\f{1}{2}})(\ln{N})^{\f{1}{2}}),
\label{Pineq}
\end{equation}
we can combine \eqref{ieq:estF} and \eqref{Pineq} to find

\begin{equation}\label{ieq:s}
s |u|_{A_{\a-1, \b}^r} < |u|_{A_{\a, \b}^r}
\end{equation}
and 

\begin{equation}\label{ieq:ext}
  \int_0^{\f{s^2}{2}}(\p_x^r u(x))^2 x^{\a+r}\e^{-\b x}\dd{x}
  < \int_{\f{s^2}{2}}^{+\infty}(\p_x^r u(x))^2 x^{\a+r}\e^{-\b x}\dd{x}.
\end{equation}
In other words, as the frequency indicator increases, the norm of
$\p_x^ru(x) \cdot \mathbb{I}_{(s^2/2,+\infty)}(x)$ becomes larger than
that of $\p_x^ru(x)\cdot \mathbb{I}_{(0,s^2/2)}(x)$, implying scaling
is indeed needed to enlarge the computational domain because $\|
\p_x^ru\cdot \mathbb{I}_{(x>s^2/2)}\|_{\o_{\a,\b}}$ is the dominant component of
$\| \p_x^ru\|_{\o_{\a,\b}}$. Here, $\mathbb{I}_S(x)$ denotes the
characteristic function on a set $S$. The verification of inequality
\eqref{ieq:ext} can be finished by contradiction.  If \eqref{ieq:ext}
does not hold, we would have

\begin{align*}
|u|_{A^r_{\a,\b}}^2 &= \int_0^{+\infty} (\p_x^r u(x))^2 x^{\a+r}\e^{-\b x}\dd{x} \\
&\le 2 \int_0^{\f{s^2}{2}}(\p_x^r u(x))^2 x^{\a+r}\e^{-\b x}\dd{x} \\
&\le 2 \cdot \frac{s^2}{2}\int_0^{\f{s^2}{2}}(\p_x^r u(x))^2 x^{\a+r-1}\e^{-\b x}\dd{x} \\
&\le s^2 \int_0^{+\infty}(\p_x^r u(x))^2 x^{\a+r-1}\e^{-\b x}\dd{x} 
= s^2 |u|_{A^r_{\a-1,\b}}^2,
\end{align*}
which would contradict the inequality \eqref{ieq:s}. Intuitively,
basis functions of higher degree decay more slowly than those of lower
degree, so an increase in the frequency indicator implies slower decay
at infinity. This slower spatial decay as time increases requires
using a larger computational domain which is achieved by decreasing $\beta$. In practice, we can also obtain good numerical results
using $\a=0$ although no theoretical result like the above observation
is guaranteed since $A_{-1, \b}^r$ is not defined.

%\section{Exterior-error-dependent moving method}
\section{Exterior-error-dependent moving}
\label{sec:moving}

Dynamics in unbounded domains can be much richer than the simple
diffusive behavior successfully captured by our frequency-dependent
scaling.  Other physical mechanisms may induce, for example,
translations (Examples~\ref{ex:mfd} and \ref{ex:nld}) and emerging
oscillations (Example~\ref{ex:refine}). A purely scaling approach
fails in these cases.

In this section, we develop an exterior-error-dependent moving method that will
be able to resolve a solution's decay in an undetermined exterior
domain $\Lambda_{\rm e}\coloneqq (x_L,+\infty)$.  Alg.~\ref{algmoving}
presents the pseudo-code of our exterior-error-dependent moving
technique. In the algorithm, we first need to determine the
time-dependent left-end point $x_L$. Next, we move the spectral basis
accordingly so that the spectral approximation for an unknown function
$u(x)$ in $\Lambda_{\rm e}$ (denoted by $U_{N,x_L}^{(\a, \b)}(x)$)
maintains accuracy. To implement this procedure, we adopt an {\em
  exterior-error indicator}:

 \begin{equation}
   \mathcal{E}(U_{N,x_L}^{(\a, \b)},x_R) = \f{\|\p_x U_{N,x_L}^{(\a, \b)}
\cdot\mathbb{I}_{(x_R,+\infty)}\|_{\o_{\a, \b}}}{\|\p_x U_{N,x_L}^{(\a, \b)}\cdot
     \mathbb{I}_{(x_L,+\infty)}\|_{\o_{\a, \b}}},
\label{errorindicator}
 \end{equation}
 which measures the proportion of the norm  $\|\p_x U_{N,x_L}^{(\a, \b)}\cdot
     \mathbb{I}_{(x_L,+\infty)}\|_{\o_{\a, \b}}$ inside a prescribed unbounded domain
  $(x_R,+\infty)$.

The subroutine $\Call{exterior\_error\_indicator}{}$ in
Lines~\ref{m:e0}, \ref{m:e1}, and \ref{m:e2} of Alg.~\ref{algmoving}
calculates $\mathcal{E}(U_{N,x_L}^{(\a, \b)},x_R)$.  Here, following
the often-used $\frac{2}{3}$-rule
\cite{Hou2007Computing,Orszag1971On}, we choose
$x_R=x_{[\f{N+2}{3}]}^{(\a, \b)}$ from the collocation points
$x_j^{(\a, \b)} (j=0,1,\ldots,N)$ in the exterior domain $\Lambda_{\rm
  e}$.

Intuitively, if $u(x)$ moves rightward in time, such as the moving
Fermi-Dirac distribution in Example~\ref{ex:mfd}, the spectral
approximation at large distances may deteriorate and the
exterior-error indicator $\mathcal{E}(U_{N,x_L}^{(\a, \b)})$ will
increase.  Consequently, the moving mechanism is triggered in
Line~\ref{m:if} of Alg.~\ref{algmoving}, and completed by updating the
left end point $x_L=x_L+d_0$ in Line~\ref{m:xL}.  Thus, the starting
point of the spectral approximation also moves rightward with time to
capture the translation.

 % of $u(x)$.

The displacement $d_0=\min\{n\delta, d_{\rm max}\}$ is determined by
the $\Call{move}{}$ subroutine in Line~\ref{m:d0}, where $n$ is the
smallest integer satisfying $\mathcal{E}(U_{N,x_L}^{(\a,
  \b)},x_R+n\delta) <\mu e_0$, $\delta$ is the minimum displacement,
$d_{\rm max}$ is the maximum displacement, and $\mu$ represents the
threshold of the increase in the exterior-error indicator that we can
tolerate. In practice, $d_{\rm max}$ should be based on a prior
knowledge of the maximum translation speed of the function $u(x)$. We usually
choose $\mu \gtrsim 1$ to prevent the exterior-error indicator from
becoming too large without invoking moving. The $\Call{move}{}$
subroutine also generates $U_{N,x_L+d_0}^{(\a, \b)}$ from
$U_{N,x_L}^{(\a, \b)}$.

\begin{algorithm}[t]
  \caption{\small Pseudo-code of spectral methods with
    exterior-error-dependent moving.}
\begin{algorithmic}[1]
\State Initialize $N$, $\Delta t$, $T$, $\alpha$, $\beta$, $U_{N,0}^{(\alpha,\beta)}(0)$, $\mu>1$, $d_{\rm max}>\delta>0$
\State $t\gets 0$
\State $x_L \gets 0$
\State $x_R \gets x^{(\alpha, \beta)}_{[\f{N+2}{3}]}$
\State $e_0 \gets \Call{exterior\_error\_indicator}{U_{N,x_L}^{(\a, \b)}(0), x_R}$\label{m:e0}
\While{$t<T$}
\State $U_{N,x_L}^{(\alpha,\beta)}(t+\Delta t)\gets \Call{evolve}{U_{N,x_L}^{(\alpha,\beta)}(t),\Delta t}$
\State $e\gets \Call{exterior\_error\_indicator}{{U_{N,x_L}^{(\alpha,\beta)}(t+\Delta{t})}, x_R}$\label{m:e1}
\If{$e > \mu e_0$}\label{m:if}
\State $(d_0, U_{N,x_L+d_0}^{(\alpha,\beta)})\gets \Call{move}{U_{N,x_L}^{(\alpha,\beta)}(t+\Delta t), \delta, d_{\rm max}, \mu e_0}$\label{m:d0}
\State $x_L \gets x_L+d_0$\label{m:xL}
\State $x_R \gets x_R+d_0$\label{m:xR}
%\State $x_L(t+\Delta{t}) \gets x_R(t)+d_0$
\State $e_0 \gets \Call{exterior\_error\_indicator}{{U_{N,x_L}^{(\alpha,\beta)}(t+\Delta{t})}, x_R}$\label{m:e2}
\EndIf
\State $t\gets t+\Delta t$
\EndWhile
\end{algorithmic}\label{algmoving}
\end{algorithm}

\begin{figure}[htp]
\begin{center}
   \includegraphics[width=6.2in]{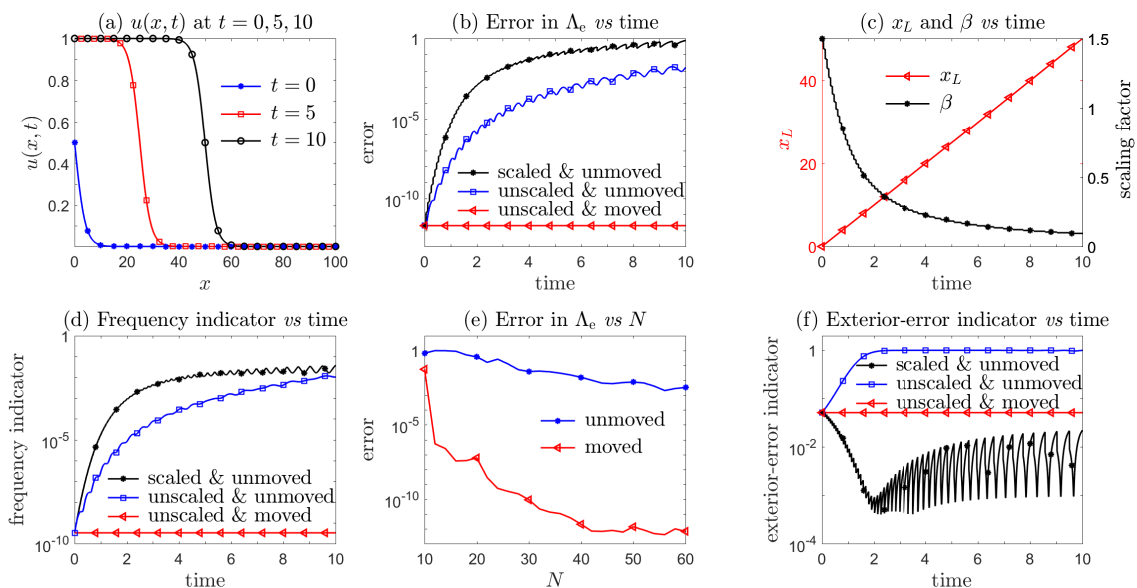}
\end{center}
\vspace{-2mm}
 \caption{\small Numerical approximation to the moving Fermi-Dirac
   distribution $u(x,t)$ given by Eq.~\eqref{eq:mfd}. The moving
   algorithm~\ref{algmoving} produces much more accurate solutions and
   recovers a faster spectral convergence with respect to the expansion
   order $N$ in the exterior domain $\Lambda_{\rm e}=(x_L, +\infty)$,
   whereas a pure scaling fails to capture this translation. The data
   in the last plot are measured at $t=10$.}
     \label{fig2}
\end{figure}

\begin{example}\label{ex:mfd}
\rm

In this example, we consider the moving Fermi-Dirac distribution
\begin{equation}\label{eq:mfd}
%u_3 = \f{1}{1+\e^{x-9-0.2t}},
u(x, t) = \f{1}{1+\e^{\f{x-5t}{2}}},
\end{equation}
which travels to the right at a speed of $5$ without any shape change 
(see Fig.~\ref{fig2}(a)). The scaling algorithm~\ref{algscaling},
equipped with the same parameters that worked well for the diffusive
Fermi-Dirac distribution in Example~\ref{ex:fd}, fails to capture the
translation. In fact, the errors of the scaled solutions are larger
than those of unscaled ones as shown in Fig.~\ref{fig2}(b).  It seems
that the decrease of the scaling factor (black curve with asterisks in
Fig.~\ref{fig2}(c)) cannot compensate for the increase in the
frequency indicator (black curve with asterisks in
Fig.~\ref{fig2}(d)). In other words, the scaling algorithm~\ref{algscaling} mistakes translation as diffusion and performs
excessive scaling. In contrast, the exterior-error-dependent moving
algorithm~\ref{algmoving} with $\delta=0.004$, $d_{\rm max}=0.04$ and
$\mu=1.005$ succeeds in producing a much more accurate approximation
to the moving Fermi-Dirac distribution given by Eq.~\eqref{eq:mfd} in the exterior domain
$\Lambda_{\rm e}$, with errors kept under $10^{-11}$ up to time $t=10$
(red curve with left-pointing triangles in Fig.~\ref{fig2}(b)). The
moving technique recovers a faster spectral convergence with respect to the
expansion order $N$ as shown in Fig.~\ref{fig2}(e).
%
%(where we have
%set $\delta=0.004$, $d_{\rm max}=0.04$ and $\mu=1.005$).
%

During the moving process, the exterior-error indicator
$\mathcal{E}(U_{N,x_L}^{(\a, \b)}, x_R)$ is well controlled (red curve
with left-pointing triangles in Fig.~\ref{fig2}(f)) and the left-end
point of the exterior domain closely tracks the uniform linear motion
(red curve with left-pointing triangles in
Fig.~\ref{fig2}(c)). The exterior-error indicator monotonically
increases for the unscaled and unmoved solutions (blue curve with
squares in Fig.~\ref{fig2}(f)) and oscillates rapidly for the scaled
and unmoved solutions (black curve with asterisks in
Fig.~\ref{fig2}(f)).  Moreover, the similarity between the relative
error and frequency indicator as a function of time is again confirmed
by comparing Fig.~\ref{fig2}(d) to Fig.~\ref{fig2}(b), thus providing
strong evidence for the effectiveness of using the frequency
indicator~\eqref{FREQSCAL}. Spectral convergence in $N$ is clearly
observed for the moving spectral method in Fig.~\ref{fig2}(e) while
the error decays slowly with $N$ for the unmoved spectral method.
\end{example}

\begin{figure}[htb]
  \begin{center}
    \includegraphics[width=6.2in]{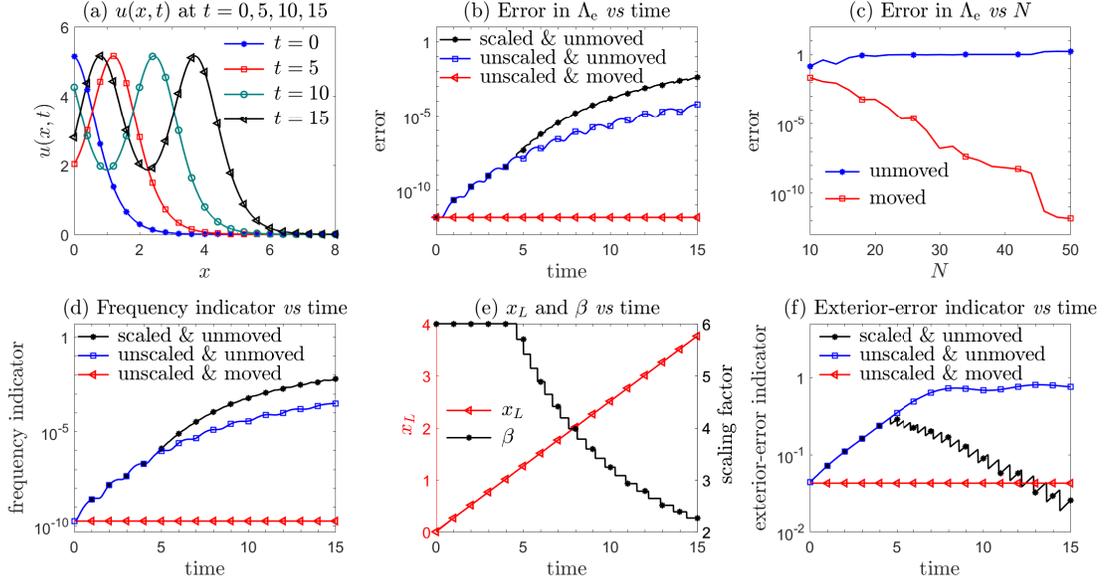}
    \end{center}
        \caption{\small Approximating a two-hump nonlinear Dirac
          solitary wave. The moving algorithm Alg.~\ref{algmoving}
          produces much more accurate solutions and recovers a faster
          spectral convergence with respect to the expansion order $N$
          in the exterior domain $\Lambda_{\rm e}=(x_L, +\infty)$,
          whereas a pure scaling approach fails to capture this
          translation. The data in the last plot are measured at
          $t=15$.}
     \label{fig3}
\end{figure}

\begin{example}\label{ex:nld}
  \rm 
  Another class of dynamical systems are described by solitons or
  solitary waves in which nonlinearities and dispersion counteract.
  While solitons have been well-studied, there has been recent
  interest in nonlinear Dirac solitary waves as they emerge naturally
  in many physical systems \cite{Cuevas2016}. Stability of the
  nonlinear Dirac solitary waves on the whole line and its connection
  to the multi-hump structure is a challenging topic of research
  \cite{Shao2014Stability,XuShaoTangWei2015,Comech2019}.  In this
  example, we approximate a right-moving two-hump solitary wave, the
  explicit form of which is given in \cite{SHAO2005interaction} with
  $v=0.25$, $\lambda=0.5$, $m=1$, $x_0=-1.5$ and $\Lambda=0.1$.
The reference solutions are plotted in Fig.~\ref{fig3}(a).

Numerical results are displayed in Fig.~\ref{fig3} where we set
$\delta=0.004$, $d_{\rm max}=0.012$, $\mu=1.005$.  It can be readily
observed there that the exterior-error-dependent moving
algorithm~\ref{algmoving} produces much more accurate solutions with
errors kept under $10^{-11}$ until the final time $t=15$ (red curve
with left-pointing triangles in Fig.~\ref{fig3}(b)). The moving
algorithm also recovers a faster spectral convergence with respect to the
expansion order $N$ (see Fig.~\ref{fig3}(c)).  The scaling-only
algorithm~\ref{algscaling} fails to maintain the accuracy (black curve
with asterisks in Fig.~\ref{fig3}(b)). The similarity between the
relative error and frequency indicator is again confirmed by comparing
Fig.~\ref{fig3}(d) to Fig.~\ref{fig3}(b).
\end{example}

In Examples~\ref{ex:mfd} and \ref{ex:nld}, the exterior-error
indicator~\eqref{errorindicator} efficiently guides us in finding an
$x_L$ such that the moved spectral approximation retains accuracy in
the resulting exterior domain.  The accuracy arises from the fact that
the exterior-error indicator is related to the upper bound of the
error for asymptotically large $x$.  If we assume a large indicator
$\mathcal{E}(U_{N,x_L}^{(\a, \b)},x_R)>\mu$ with $\mu\in(0,1)$, then
the upper bound for the error in $x>x_{R}$ is larger than the upper
bound for the error in $\Lambda_{\rm e}$:

\begin{align*}
\mathcal{E}(U_{N,x_L}^{(\a, \b)},x_R)>\mu &\Rightarrow |U_{N,x_L}^{(\a, \b)} \cdot \mathbb{I}_{[x_R,+\infty)}|_{A_{\a-1, \b}^1}>\mu |U_{N,x_L}^{(\a, \b)}
\cdot\mathbb{I}_{[x_L,+\infty)}|_{A_{\a-1, \b}^1}, \\
& \Rightarrow
 |U_{N,x_L}^{(\a, \b)}\cdot\mathbb{I}_{[x_R,+\infty)}|_{A_{\a, \b}^1}
>\mu |U_{N,x_L}^{(\a, \b)} \cdot \mathbb{I}_{[x_L,+\infty)}|_{A_{\a, \b}^1}.
\end{align*}
The solution in the interior domain $\Lambda_{\rm i}\coloneqq (0,x_L]$
is not approximated by the basis functions used to approximate the
solution in the exterior domain.  Obstacles to designing moving mesh
methods in unbounded domains include the construction of an {\em
  interior numerical solution} and its consistent coupling with the
{\em exterior spectral approximation}. More on these issues will be
illustrated in Example~\ref{ex:refine}.

\begin{figure}[htb]
  %\vspace{-2cm}
    \begin{center}
    \includegraphics[width=4.4in]{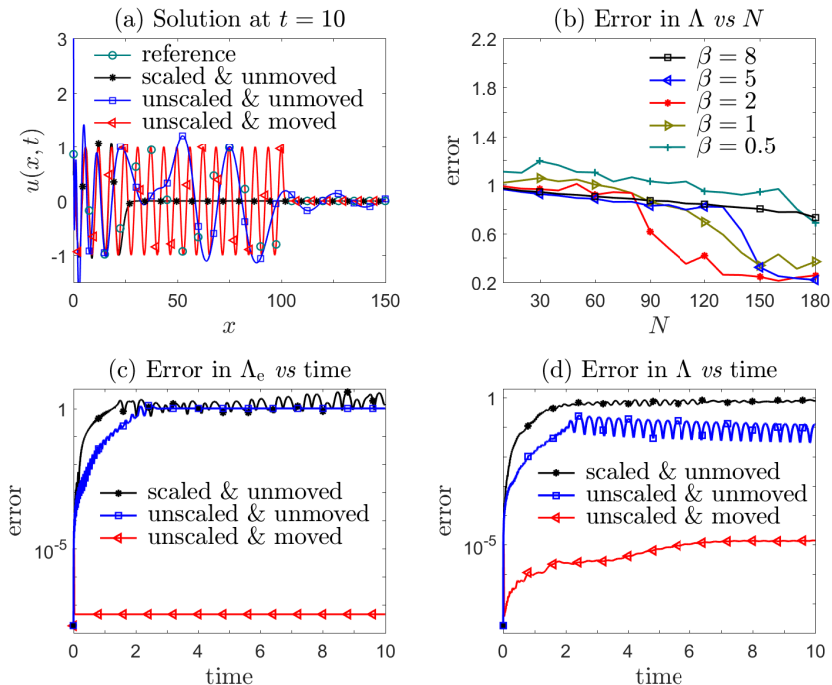}
    \end{center}
        \caption{\small Oscillations emanate from the left but the
          moving algorithm~\ref{algmoving} generates accurate
          solutions in the exterior domain $\Lambda_{\rm e}$, with
          relative errors under $10^{-7}$ up to $t=10$ with $N=30$
          (red curve with left-pointing triangles in (c)).  By further
          coupling with a spectral approximation using $80$ Chebyshev
          polynomials in the interior domain $\Lambda_{\rm i}$, we
          generate the whole solution with total relative error, up
          until $t=10$, under $2\times 10^{-5}$, as shown by the red
          curves with left-pointing triangles in (a) and (d). The data
          in (b) are measured at $t=10$.}
     \label{fig4}
\end{figure}

\begin{example}\label{ex:refine}
\rm

Let us approximate the following function in $\Lambda$: 
\begin{equation}\label{eq:refine}
u(x, t) =
\left\{
\begin{aligned}
\cos(x-10t), & \quad  x \leq 10t,\\
\e^{-(x-10t)^2}, & \quad  x > 10t,
\end{aligned}
\right.
\end{equation}
which represents a wave with period $2\pi$ traveling to the right with
speed 10 and exponentially decaying at infinity. The reference
solution $u(x,10)$ is plotted by the green curve with circles in
Fig.~\ref{fig4}(a), which coincides with the red curve with
left-pointing triangles hat approximates $u$ seperately in $\Lambda_i$ 
and $\Lambda_e$ using different basis funcitons.  As shown by the blue
curve with squares in Fig.~\ref{fig4}(a), applying a Laguerre spectral
approximation with $N=30$ and $\beta=5$ in $\Lambda$ fails to
accurately approximate $u(x,t)$.  This failure arises because more
oscillations emerge from $x=0$ and translate to $+\infty$ as time
evolves. Specifically, at $t=10$, the reference solution $u(x,t)$
possesses $32$ extrema while any Laguerre spectral
approximation~\eqref{SpecD} with $N=30$ can have at most $30$ extrema,
implying that the approximation is doomed to fail since all
oscillations cannot be captured. Simply increasing the number of basis
functions does little to help, even with different scaling factors as
shown in Fig.~\ref{fig4}(b). The ineffectiveness of increasing $N$ is
mainly due to the presence of oscillatory components with
significantly different frequencies in each of the two different
domains. As shown by the black curves with asterisks in
Figs.~\ref{fig4}(a, c, d), the scaling
technique is also doomed to fail because it totally neglects this
scale difference and only adjusts the scaling factor to redistribute
collocation points.
%
%the size of which is so fixed that the resulting spectral
%approximation cannot capture all emerging oscillations from the left.
%

We propose a divide-and-conquer strategy to address
Example~\ref{ex:refine} that can be implemented by applying two
subroutines, within each time step.  The first step is to use the
exterior-error-dependent moving algorithm~\ref{algmoving} to determine
the exterior spectral approximation for the exponential decay
component of the reference solution.  The second step is to introduce
a new spectral approximation in the remaining bounded interior domain
$\Lambda_{\rm i}$ for the left-side oscillating component. The full
numerical solution in the half-line $\Lambda$ is constructed from
concatenating the solution in the exterior domain $\Lambda_{\rm e}$ to
the one in the interior domain $\Lambda_{\rm i}$.

Fig.~\ref{fig4}(c) plots the error in the exterior domain against time
and shows that the errors of of the moved solution with $N=30$,
$\delta=0.008$, $d_{\rm max}=0.08$ and $\mu=1.001$ are kept under
$10^{-7}$ up to time $t=10$ (red curve with left-pointing triangles),
confirming that the Laguerre spectral approximation is accurate in the
exterior domain.  In fact, the numerical values of $x_L$ obtained by
the moving algorithm~\ref{algmoving} are consistent with the expected
value of $10t$ as shown in Eq.~\eqref{eq:refine}.
%
%, for example, $x_l=100.032$ for $t=10$.
%
Coupling the exterior solution with a spectral approximation using
$80$ Chebyshev polynomials in the interior domain, we find a combined
numerical solution with total relative error under $2\times 10^{-5}$
up to $t=10$ (red curves with left-pointing triangles in
Figs.~\ref{fig4}(a, d)) using $111=31+80$ total basis
functions.  By contrast, Fig.~\ref{fig4}(b) shows that the errors for
direct refinement using $N=180$ are larger than $0.2$.

It must be pointed out that when solving PDEs in unbounded domains, we
may need information about the solution in the exterior domain to
construct the interior numerical solution. Further discussion on this
point can be found in Example~\ref{ex:PDEsolve}.
\end{example}

%%%%%%%%%%%%%%%%%%%%%%%
%%%%%%%%%%%%%%%%%%%%%%%%%%%

\section{Spectral methods incorporating both scaling and moving}
\label{sec:scalingmoving}

For problems that involve both translation and diffusion in unbounded
domains, we need to incorporate both the moving and scaling
procedures.  Since the scaling algorithm~\ref{algscaling} may
mistake translation for diffusion and trigger an inappropriate scaling
as shown in Examples \ref{ex:mfd} and \ref{ex:nld}, we propose a
``first moving then scaling'' algorithm.  The associated pseudo-code
is described in Alg.~\ref{algmovingscaling}. A direct application of
Alg.~\ref{algmovingscaling} to Example~\ref{ex:fd} recovers exactly
the same results as Alg.~\ref{algscaling} since the moving procedure
is not invoked.  When Alg.~\ref{algmovingscaling} is applied to
Examples \ref{ex:mfd} and \ref{ex:nld}, it gives the same results as
Alg.~\ref{algmoving} since the scaling mechanism is not
triggered. That is, the combined moving-scaling algorithm~\ref{algmovingscaling} can deal with both translation-only and
diffusion-only problems since it can distinguish translation from
diffusion.

Alg.~\ref{algmovingscaling} can be extended to
unbounded domains in multiple dimensions in a dimension-by-dimension
manner by using the tensor product of one-dimensional basis
functions. For example, consider the two-dimensional spectral approximation
\begin{equation}
  U_{N, x_L, y_L}^{(\vec{\alpha}, \vec{\beta})}(x, y)\coloneqq
  \sum\limits_{\ell=0}^{N_x}\sum\limits_{m=0}^{N_y}
  u^{(\vec{\alpha}, \vec{\beta})}_{\ell, m}
  \mathcal{L}_{\ell}^{\alpha_x, \beta_x}(x)\mathcal{L}_{m}^{\alpha_y, \beta_y}(y)
\end{equation}
in $\Lambda_{\rm e}^x\times\Lambda_{\rm e}^y: =(x_L, +\infty)\times(y_L,
+\infty)$ where $\vec{\alpha}=(\alpha_x,\alpha_y)$ and
$\vec{\beta}=(\beta_x,\beta_y)$.  We choose the exterior-error
indicator in $x$-dimension to be
\begin{align}
  \mathcal{E}_x(U_{N, x_L, y_L}^{(\vec{\alpha}, \vec{\beta})}(x, y), x_R)& \coloneqq
  \mathcal{E}(\tilde{U}_{N, x_L}^{(\alpha_x, \beta_x)}(x), x_R), \\
  \tilde{U}_{N, x_L}^{(\alpha_x, \beta_x)}(x) & \coloneqq
  \int_{\Lambda_{\rm e}^y} U_{N, x_L, y_L}^{(\vec{\alpha}, \vec{\beta})}(x, y) \dd{y}.
\end{align}
Similarly, $\mathcal{E}_y(U_{N, x_L, y_L}^{(\vec{\alpha},
  \vec{\beta})}(x, y), y_R)$ gives the exterior-error indicator in
$y$-dimension. Accordingly, we use $\mathcal{E}_x(U_{N, x_L,
  y_L}^{(\vec{\alpha}, \vec{\beta})}, x_R)$ to judge the \t{if}
statement in Line \ref{move:if} of Alg.~\ref{algmovingscaling}.  If
satisfied, then the $\Call{move}{}$ subroutine in Line \ref{move} will
move the solution in $x$-direction via $x_L \to x_L
+d_{0}^x$. Simultaneously, we use $\mathcal{E}_y(U_{N, x_L,
  y_L}^{(\vec{\alpha}, \vec{\beta})}(x, y), y_R)$ to determine the
shift in the $y$-direction.

To allow scaling in $x$-direction, the corresponding frequency
indicator can be defined as
\begin{equation}
  \mathcal{F}_x(U_{N, x_L, y_L}^{(\vec{\alpha}, \vec{\beta})})\coloneqq
  \left(\f{\sum\limits_{\ell=N_x-M_x+1}^{N_x}\sum\limits_{m=0}^{N_y}\g_{\ell}^{(\a_x, \b_x)}
    \g_m^{(\a_y, \b_y)} (u^{(\vec{\alpha}, \vec{\beta})}_{\ell, m})^2}
       {\sum\limits_{\ell=0}^{N_x}\sum\limits_{m=0}^{N_y}
         \g_{\ell}^{(\alpha_x, \beta_x)}
         \g_m^{(\a_y, \b_y)}(u^{(\vec{\alpha}, \vec{\beta})}_{\ell, m})^2}\right)^{\f{1}{2}},
\label{freqx}
\end{equation}
where $M_x = [\f{N_x}{3}]$ and $N_x, N_{y}$ are the expansion orders in the $x$-, $y$-directions, respectively.  Similarly, we define
$\mathcal{F}_y$ to be the frequency indicator in $y$-direction.  We
first keep $\beta_y$ fixed and use $\mathcal{F}_x$ to evaluate the
\textbf{if} statement in Line \ref{scale:if} for scaling. If scaling
in $x$-direction is needed, then the \t{while} loop in Line
\ref{scale:while} will update the scaling factor to $\tilde{\beta}_x$.
Simultaneously, we fix $\beta_x$ and use $\mathcal{F}_y$ to update the
scaling factor in the $y$-direction to $\tilde{\beta}_y$. After that,
the scaling factors for time $t+\Delta{t}$ are set to
$\tilde{\beta}_x$ and $\tilde{\beta}_y$.

\begin{figure}[htb]
 %\vspace{-2cm}
    \begin{center}
    \includegraphics[width=4.4in]{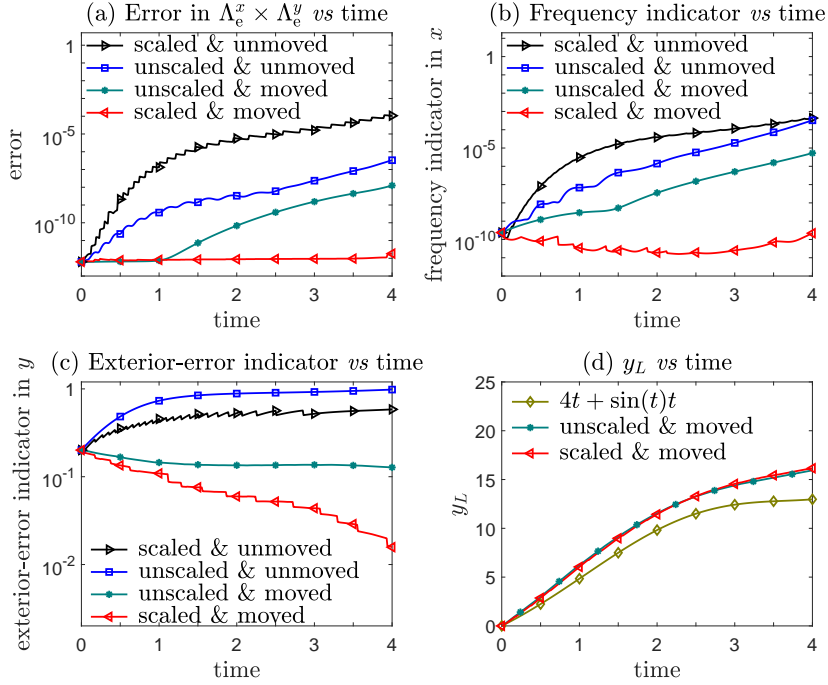}
	\end{center}
        \caption{\small A two-dimensional oscillatory function with
          both translation and diffusion given by
          Eq.~\eqref{ex:2d1}. Only the combined moving-scaling
          algorithm~\ref{algmovingscaling} produces accurate solutions
          in the exterior domain with errors kept under $10^{-11}$ up
          to $t=4$. The need for
          combining moving and scaling is evident. For simplicity, we
          only used $\mathcal{F}_x$ (the frequency indicator in the $x$-direction),
          $\mathcal{E}_y$ (the exterior-error indicator in the
          $y$-direction), and $y_L$ (the left end of
          $\Lambda_{\rm e}^y$) as an example.  The corresponding
          curves for $\mathcal{F}_y$, $\mathcal{E}_x$, and $x_L$ are
          very similar and not shown.  Here, we used $N_x=N_y=40$,
          and the initial scaling factors: $\beta_x=\beta_y=2.5$.}
     \label{fig5}
\end{figure}

\begin{algorithm}
\caption{\small Pseudo-code of spectral methods with both scaling and moving.}
\begin{algorithmic}[1]
\State Initialize $N$, $\nu>1$, $q<1$, $\Delta t$, $T$, $\alpha$, $\beta$, $U_{N}^{(\alpha,\beta)}(0)$, $\underline{\b}$, $\mu>1$, $d_{\rm max}>\delta>0$,  $x_R(0) = x^{(\alpha, \beta)}_{[\f{N+2}{3}]}$
\State $t, x_L \gets 0$
\State $x_R \gets x^{(\alpha, \beta)}_{[\f{N+2}{3}]}$
\State $f_0 \gets \Call{frequency\_indicator}{U_{N, x_L}^{(\alpha,\beta)}(x, t)}$
\State $e_0 \gets \Call{exterior\_error\_indicator}{U_{N, x_L}^{(\a, \b)}(0), x_R}$
\While{$t<T$}
\State $x_R \gets x^{(\alpha, \beta)}_{[\f{N+2}{3}]}$\label{update_xR}
\State $U_{N, x_L}^{(\alpha,\beta)}(x, t+\Delta t)\gets \Call{evolve}{U_{N,x_L}^{(\alpha,\beta)}(x, t)),\Delta t}$\label{intermediate}
\State $e\gets \Call{exterior\_error\_indicator}{{U_{N, x_L}^{(\alpha,\beta)}(x, t+\Delta{t})},x_R}$
\If{$e > \mu e_0$} \label{move:if}
\State $(d_0, U_{N, x_L+d_0}^{(\a, \b)})\gets \Call{move}{U_{N, x_L}^{(\alpha,\beta)}(x, t+\Delta t), \delta, d_{\rm max}, \mu e_0}$ \label{move}
\State $x_L \gets x_L+d_0$  \label{move:d0}
\State $e_0 \gets \Call{exterior\_error\_indicator}{{U_{N, x_L}^{(\alpha,\beta)}(x, t+\Delta{t})},x_R}$
\EndIf
%\State $e_0 \gets e$
\State $f\gets \Call{frequency\_indicator}{U_{N, x_L}^{(\alpha,\beta)}(x, t+\Delta t)}$
\If{$f > \nu f_0$}\label{scale:if}
\State $\tilde{\b}  \gets q \beta$
\State $U_{N, x_L}^{(\a,\tilde{\b})}  \gets  \Call{scale}{U_{N, x_L}^{(\alpha,\beta)}(x, t+\Delta t),\tilde{\b}}$
\State $\tilde{f} \gets \Call{frequency\_indicator}{U_{N, x_L}^{(\a,\tilde{\b})}}$
\While{$\tilde{f}\leq f$ and $\tilde{\b}\geq \underline{\b}$}\label{scale:while}
\State $\b \gets \tilde{\b}$
\State $U_{N, x_L}^{(\alpha,\beta)}(x, t+\Delta t)\gets U_{N, x_L}^{(\a,\tilde{\b})}$
\State $f_0 \gets \tilde{f}$
\State $f \gets \tilde{f}$
\State $\tilde{\b}  \gets q \beta$
\State $U_{N, x_L}^{(\a,\tilde{\b})}  \gets  \Call{scale}{U_{N, x_L}^{(\alpha,\beta)}(x, t+\Delta t),\tilde{\b}}$
\State $\tilde{f} \gets \Call{frequency\_indicator}{U_{N, x_L}^{(\a,\tilde{\b})}}$
\EndWhile
\EndIf
\State $t\gets t+\Delta t$
\EndWhile
\end{algorithmic}\label{algmovingscaling}
\end{algorithm}

\begin{example}\label{ex:2D}
\rm
We will investigate the performance of Alg.~\ref{algmovingscaling}
in a two-dimensional unbounded domain by considering the function 
\begin{equation}\label{ex:2d1}
  u(x, y, t) = \cos(\f{xy}{400})\cdot
  \f{1}{1+\e^{\f{x-6t-2-t\cos(t)}{2+0.3t}}}\cdot
  \f{1}{1+\e^{\f{y-4t-2-t\sin(t)}{2+0.4t}}}, \,\,\, x, y, t>0,
\end{equation}
which displays both advective and diffusive behavior.  This function
exhibits oscillations in space from the factor $\cos(\f{xy}{400})$, an
exponential decay, and a translation to infinity with time-varying
velocity $\vec{v}=(v_x,v_y) = (6+\cos(t),4+\sin(t))$.  The numerical
results shown in Fig.~\ref{fig5} are generated using a time step
$\Delta{t}=0.01$, the same parameters in the $x$-, $y$- directions,
and $N_x=40$, $\mu_x=1.003$, $\delta_x= 0.005$, $d_{\rm max}^x=0.1$.

As expected, only the combined scaling-moving algorithm~\ref{algmovingscaling} keeps the errors in the exterior domain
under $10^{-11}$ (up to the final time $t=4$), as shown by the error
curves in Fig.~\ref{fig5}(a). This accuracy is achieved because the
corresponding frequency indicator and exterior-error indicator are
controlled by our ``first moving then scaling'' techniques, see e.g.,
$\mathcal{F}_x$ in Fig.~\ref{fig5}(b) and $\mathcal{E}_y$ in
Fig.~\ref{fig5}(c).

Although the moving algorithm~\ref{algmoving} may accurately
capture the function near the left end of the exterior domain, 
the resulting exterior-error indicator does not stay low enough to
preserve accuracy in the exterior domain $\Lambda_{\rm
  e}^x\times\Lambda_{\rm e}^y$, as shown by the green curves with
asterisks in Figs.~\ref{fig5}(a, c, d). 
The moving algorithm neglects the diffusion and thus uses an improper
(smaller) $x_R$ and $ y_R$. The right choice for these two variables
depends on proper scaling for the diffusion, revealing why we need to
update $x_R$ in Line \ref{update_xR} of Alg.~\ref{algmovingscaling}
after scaling. That is, the moving determines $x_L$ while the scaling
determines $x_R$, making it necessary to combine moving with scaling.

\end{example}

As we have mentioned in Example~\ref{ex:refine}, numerically solving
evolving PDEs in unbounded domains requires both the interior solution
$U_{x_L(t)}^{\text{interior}}(x, t)$ in $\Lambda_{\rm
  i}(t)=(0,x_L(t)]$ and the exterior solution $U_{N, x_L(t)}^{(\a,
   \b)}(x, t)$ in $\Lambda_{\rm e}(t)=(x_L(t),+\infty)$ after applying
 the divide-and-conquer strategy. When using the moving-scaling
 algorithm~\ref{algmovingscaling} to march the solution from $t$
 to $t+\Delta t$, if the moving mechanism is not triggered
 (\textit{i.e.}, $x_L$ is unchanged), then the interior and exterior
 solutions can be updated individually in the normal way. If it is
 triggered, extra steps are needed to approximate the solution in the
 enlarged interior domain $\Lambda_{\rm i}(t+\Delta{t})=\Lambda_{\rm
   i}(t)\cup(\Lambda_{\rm e}(t)\setminus\Lambda_{\rm e}(t+\Delta{t}))$
 since $x_L(t+\Delta t) = x_L(t)+d_0$ after running Line~\ref{move:d0}
 of Alg.~\ref{algmovingscaling}.

 In the next Example, we will test the ability of
 Alg.~\ref{algmovingscaling} to solve a one-dimensional PDE where we
 will use the intermediate (unmoved) exterior solution $U_{N,
   x_L(t)}^{(\a, \b)}(x, t+\Delta t)$ (obtained immediately after
 running Line~\ref{intermediate}) to interpolate the required function
 values in $\Lambda_{\rm i}(t+\Delta{t}) \setminus \Lambda_{\rm
   i}(t)$.

\begin{example}\label{ex:PDEsolve}
\rm We solve the following first-order PDE
\begin{equation}\label{PDEMOVESCALE}
\partial_t u(x,t) + \left(2 + \f{x-2t}{2+t}\right)\partial_x u(x,t) = 0
\end{equation}
with initial data $u(x, 0) = ({1+\mbox{e}^{\f{x}{2}}})^{-1}$ and
Dirichlet boundary condition $u(0, t) = ({1+\mbox{e}^{\f{-2t}{2+t}}})^{-1}$.
The analytical solution is a moving and diffusive Fermi-Dirac
distribution: $u(x, t) = ({1+\mbox{e}^{\f{x-2t}{2+t}}})^{-1}$, which
travels rightward to infinity at a speed of $2$.  A simple numerical
scheme for evolving Eq.~\eqref{PDEMOVESCALE} is employed here for
testing the performance of Alg.~\ref{algmovingscaling} within the
divide-and-conquer strategy.

Specifically, we adopt the Laguerre spectral approximation
\eqref{SpecD} in the exterior domain, the first order backward finite
difference method in the interior domain, and the second order
improved Euler scheme in time. We use a nonuniform mesh, e.g., $10$
Gauss-Lobatto points, to avoid possible poor resolution in the tiny
interior domain $0<x_L<d_{\rm max}$ at short times. For $x_L\ge d_{\rm
  max}$, a uniform mesh with spacing $\Delta x=\delta=0.02$ is used so
new grid points in $\Lambda_{\rm i}(t+\Delta{t}) \setminus
\Lambda_{\rm i}(t)$ can be easily added. The other parameters were set
to $N=40$, $\mu=1.004$, $d_{\rm max}=0.2$, and $\Delta{t}=0.001$.

The results summarized in Fig.~\ref{fig6} clearly show that, up to the
final time $t=5$, the proposed divide-and-conquer strategy maintains
the errors in the whole domain $\Lambda=\Lambda_{\rm
  i}\cup\Lambda_{\rm e}$ under $2\times10^{-4}$ (red curve with
left-pointing triangles in Fig.~\ref{fig6}(a)). 
Alg.~\ref{algmovingscaling} succeeds in capturing the translation, as
shown by the red curve with left-pointing triangles in
Fig.~\ref{fig6}(b), thus determining the exterior domain $\Lambda_{\rm
  e}$.  Without this strategy, a straightforward use of the Laguerre
spectral approximation in $\Lambda$ leads to huge errors as indicated
by the blue curve with right-pointing triangles in
Fig.~\ref{fig6}(a).

Fig.~\ref{fig6}(c) shows that the frequency indicator is always kept
under $3\times10^{-10}$  as shown by the black curve with asterisks, a
sufficiently small lower error bound for scaling, by continually
shrinking the scaling factor shown as the black curve with asterisks
in Fig.~\ref{fig6}(b). The exterior-error indicator is always
maintained around $0.2$ as shown by the red curve with left-pointing
triangles in Fig.~\ref{fig6}(c), which implies the error in $(x_R,
+\infty)$ divided by the error in $\Lambda_{\rm e}$ is almost
  unchanged, ensuring small errors at infinity.  Fig.~\ref{fig6}(d)
plots $|U(x,t)-u(x,t)|$ at different times ($U(x,t)$ and $u(x,t)$
denote the numerical and analytical solution, respectively).  There is
a clear divide near $x_L$ arising from the different numerical
treatments between the interior and exterior domains.

% by the `first moving then scaling' technique built in .  the
%moving-scaling algorithm~\ref{algmovingscaling} maintains the errors
%in the whole domain under $10^{-4}$ until the final time $t=5$,

 \begin{figure}[t]
 %\vspace{-2cm}
      \begin{center}
    \includegraphics[width=4.4in]{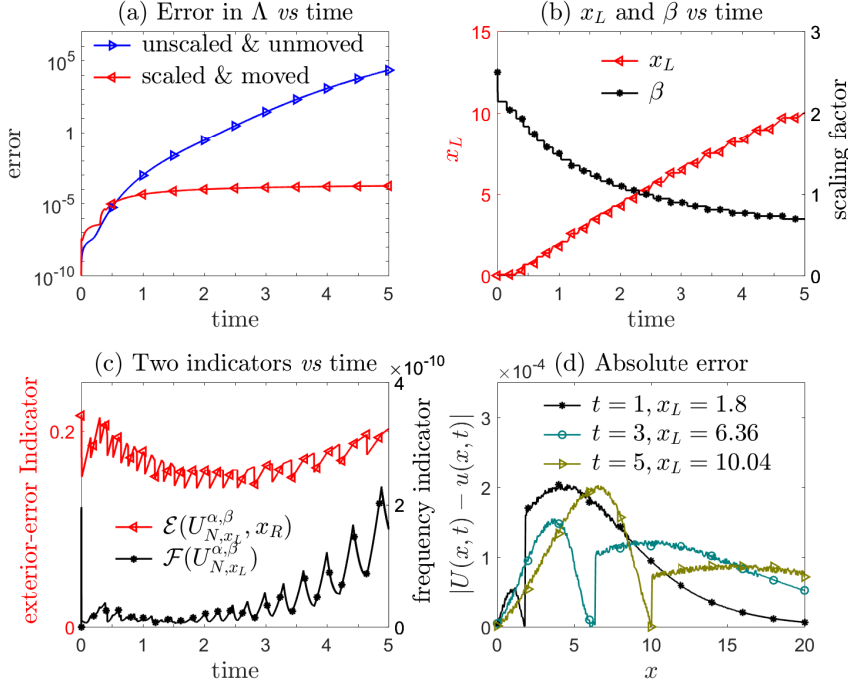}
	\end{center}
        \caption{\small Numerical results obtained by the
          moving-scaling algorithm~\ref{algmovingscaling} for the
          one-dimensional problem in Eq.~\eqref{PDEMOVESCALE}.  The
          proposed divide-and-conquer strategy maintains the errors in
          the whole domain $\Lambda=\Lambda_{\rm i}\cup\Lambda_{\rm
            e}$ under $2\times10^{-4}$ until the final time $t=5$
          where the exterior domain $\Lambda_{\rm e}$ is determined by
          the ``first moving then scaling'' technique built in to
          Alg.~\ref{algmovingscaling}. We adopt the Laguerre spectral
          approximation \eqref{SpecD} with $N=40$ in the exterior
          domain $\Lambda_{\rm e}=(x_L,+\infty)$, the first order
          backward finite difference method with spacing $\Delta
          x=0.02$ in the interior domain $\Lambda_{\rm i}=(0,x_L]$,
        and the second order improved Euler time marching scheme with
        $\Delta t=0.001$. The last plot displays the absolute
        difference between the numerical solution $U(x,t)$ and the
        analytical one $u(x,t)$ at different times.}
     \label{fig6}
\end{figure}
\end{example}

\section{Performance comparison in solving parabolic PDEs}
\label{sec:com}

We now apply the frequency-dependent scaling algorithm~\ref{algscaling} to solve 

\begin{equation}
\partial_t u(x,t) - \partial_{xx} u(x,t) = f(x,t)
\label{PARA}
\end{equation}
in $\mathbb{R}\times\Lambda$, and compare our results with those
obtained with the time-dependent scaling method developed in
\cite{MA2005}.  First, we need to generalize our scaling approach from
$\Lambda$ to $\mathbb{R}$ by using scaled Hermite polynomials, denoted
by $\m{H}_{\ell}^{(\b)}(x)$, which are mutually orthogonal under the
weight function $\o_{\b}(x)=\e^{-(\b{x})^2}$ ($\b > 0$). Similarly, we
use $\beta$ to denote the {\em scaling factor} and the frequency
indicator defined in Eq.~\eqref{FREQSCAL} still serves as a lower
bound for the interpolation error.

We use a standard Galerkin Hermite spectral method to find a solution
$U_N^{(\b)}=\sum_{\ell=0}^Nu_{\ell}^{(\b)}\hat{\mathcal{H}}_{\ell}^{\beta}(x)$
in $V_N^{(\b)} = \text{span}\{\hat{\m{H}}_0^{(\b)}(x), ...,
\hat{\m{H}}_N^{(\b)}(x)\}$ satisfying the initial condition and

\begin{equation}
(\p_tU_N^{(\b)}, v) + (\p_xU_N^{(\b)}, \p_xv)  = (f, v),\quad \forall\, v\in V^{(\b)}_N,
\label{WEAKPARA}
\end{equation}
where $\hat{\m{H}}_{\ell}^{(\b)}(x)\coloneqq
\sqrt{\o_{\b}(x)}\m{H}_{\ell}^{(\b)}({x})/\|\m{H}_{\ell}\|_{\o_{\b}}$
denotes the corresponding scaled Hermite functions and $(\cdot,
\cdot)$ is the conventional inner product in $L^2(\mathbb{R})$ space.
The Galerkin discretization~\eqref{WEAKPARA} is stable in the sense
that
\begin{equation}
  (\p_x U_N^{(b)}, \p_x U_N^{(\b)}) =
  \sum_{\ell=0}^{N+1}\f{\ell+1}{2}(u_{\ell}^{(\b)})^2-\sum_{\ell=0}^{N-2}
  \sqrt{(\ell+1)(\ell+2)}u_{\ell}^{(\b)}u_{\ell+2}^{(\b)}
\end{equation}
is strictly positive and can be controlled by
$(N+1)\|U_N^{(\b)}\|_2^2=(N+1)\sum_{\ell=0}^{N}(u_{\ell}^{(\b)})^2$.
By contrast, a time-dependent scaling factor:

\begin{equation}
\beta(t) = \f{1}{2\sqrt{\delta_0(\delta{t}+1)}}
\label{MASCALE}
\end{equation}
was taken in \cite{MA2005} to fix the instability of the
Petrov–Galerkin discretization by tuning the parameters $\delta_0$ and
$\delta$.

\begin{example}\label{ex:Hermite}
\rm 
We apply the frequency-dependent scaling
algorithm~\ref{algscaling} to Example 6.1 in \cite{MA2005}. In order
to facilitate comparison, we also adopt the same second order-accurate
Crank-Nicholson scheme to march Eq.~\eqref{WEAKPARA}, and the same
errors $E_N$ and $E_{N, \infty}$ to measure the accuracy.
Table~\ref{tab:parabolicerrorordr} presents the numerical errors with
different time steps and expansion orders where the second-order
accuracy in time and the spectral convergence in space are clearly
demonstrated. Table~\ref{tab:parabolic} compares the errors $E_N$
without scaling to those obtained using the scaling
algorithm~\ref{algscaling} and the time-dependent scaling method in
\cite{MA2005} on the same mesh.  Both scaling methods produce much
more accurate numerical results but the proposed frequency-dependent
scaling keeps the errors around or below $10^{-7}$, outperforming the
time-dependent scaling of \cite{MA2005}.

The scaling factor adjusted adaptively by the frequency indicator
\eqref{FREQSCAL} takes on the value $\beta = 0.5357$ at $t=1$ for all
choices of time steps shown in Table~\ref{tab:parabolic} whereas the
time-dependent scaling factor in \cite{MA2005} decreases to $\beta =
0.3536$ at $t=1$ (Eq.~\ref{MASCALE}). The smaller scaling factor
arises from the stability requirement $\beta'(t)+2\beta^3(t)\leq 0$,
an initial value of $0.5$, and using $\delta_0=\delta=1$ in
Eq.~\ref{MASCALE} \cite{MA2005}, and prevents the error from
decreasing when the time step is refined from $1/4000$ to $1/16000$
(see the third column of Table~\ref{tab:parabolic}).
There is no accuracy improvement without
scaling when the timestep is decreased  as shown in the second
  column of Table~\ref{tab:parabolic} where a scaling factor is fixed
  to $\beta=0.85$. 
  Regardless of what time step is used in the
unscaled method, the error $E_N$ experiences a sudden increase across
$t\in[0.3, 0.7]$, rising from below $10^{-6}$ to about $10^{-4}$, as
it fails to capture the diffusion. A similar observation was shown in
Table 6.1 of \cite{MA2005}.

\begin{table}
\centering
\caption{\small Numerical results for the parabolic problem in
  Eq.~\eqref{PARA}: Errors associated with the frequency-dependent scaling
  algorithm~\ref{algscaling} at $t=1$ with different time step
  and expansion order $N$.}
\begin{tabular}{|l|r|r|r|r|r|r|}
\hline
Time step & $N$ &$E_N(1)$ & Order & $E_{N, \infty}(1)$ & Order\\
\hline
$10^{-1}$ & \multirow{4}{*}{25} &2.500e-04 & & 2.182e-04& \\ 
$10^{-2}$ & &2.499e-07 & -2.000& 2.227e-06&1.991\\
$10^{-3}$ & &2.500e-09 & -2.000 &2.227e-08 &-2.000\\
$10^{-4}$ & &2.555e-10 & -1.991&2.350e-10 &-1.977\\
\hline
\multirow{4}{*}{$1/40000$}&$10$ & 2.203e-04& & 1.619e-04 & \\ 
&$15$& 2.189e-07&$N^{-16.85}$ & 4.335e-08 & $N^{-20.29}$\\
&$20$ &1.353e-09 & $N^{-17.68}$ & 8.880e-09 & $N^{-13.52}$\\
&$25$ &4.840e-11 & $N^{-14.93}$ & 6.183e-11 & $N^{-11.94}$\\
\hline
\end{tabular}%
\label{tab:parabolicerrorordr}%
\end{table}

\begin{table}
\centering
\caption{\small Numerical results for the parabolic problem in
  Eq.~\eqref{PARA}: Comparison of the errors at $t=1$ with
  $N=20$.}
\begin{tabular}{|p{50pt}|p{50pt}|p{73pt}|p{93pt}|}
\hline
Time step & No scaling & Time-dependent scaling in \cite{MA2005} & Frequency-dependent 
scaling 
in Alg.~\ref{algscaling}\\
\hline
1/250 & 3.969e-04 & 2.598e-06 & 3.998e-07\\ 
\hline
1/1000 & 3.910e-04 & 1.189e-06 & 2.503e-08\\
\hline
1/4000 & 3.390e-04 & 1.117e-06 & 2.085e-09\\
\hline
1/16000 & 3.390e-04 & 1.117e-06 & 1.381e-09 \\
\hline
\end{tabular}
\label{tab:parabolic}
\end{table}

\end{example}

\section{Applications to structured cell population models}
\label{sec:bio}

One example of an application requiring the solution of PDEs in an
unbounded domain is the structured population models that track
populations of cells endowed with attributes such as their size.  The
standard sizer-timer model for the density of cells with age near $a$
and size near $x$ is formulated in \cite{Metz1986}, and generalizations
to include stochasticity in growth rate is studied in 
\cite{Steinsaltz2011Derivatives, caswell2005sensitivity}. Here we address
 a continuum model describing a stochastic model for cell populations 
\cite{Xia2020c}:

\begin{equation}
\frac{\partial{n}}{\partial{t}} + \frac{\partial{n}}{\partial{a}} +
\frac{\partial{(ng)}}{\partial{x}} - \f{1}{2}\f{\p^2(\sigma{n})}{\p
  x^2} =-D(a, x, t) n(a, x, t), \quad (a,
x)\in\Lambda\times\Lambda,
\label{PDEForm}
\end{equation}
where $n(a,x,t)$ describes the density of cells with respect to age
$a$ and size $x$ at time $t$, $g(a,x,t)$ is the mean growth rate of an
individual cell and $\sigma(a,x,t)$ is the variance of stochasticity
in the growth rate, \textit{i.e.},  $\dd{x} = g\dd{t} + \sigma\dd{B}_t$, for an
individual cell.  The fluctuating growth rate manifests itself as a
diffusive term.  The right-hand-side of Eq.~\eqref{PDEForm} represents
cell division occurring with division rate $D(a,x,t)$.  Dirichlet
boundary conditions are imposed at $x=0$, $n(a, 0, t) = n_0(a, t)$,
and at $x=+\infty$, $n(a, +\infty, t)=0$ if we assume that there are no
cells of infinite size. More importantly, the boundary condition at
$a=0$ should account for two daughter cells (one of size $x$ and one
of size $y-x$) from the binary fission of a mother cell of size $y>x$:

\begin{equation}
n(x, 0, t)=2\int_0^{+\infty}\!\!\dd a\int_x^{+\infty}\!\!
\dd{y}\, \tilde{D}(a, y, x, t)n(a, y, t), 
\label{BOUNDARY}
\end{equation}
where $\tilde{D}(a, y, x, t)$ is the differential division rate
representing the rate that a cell of age $a$ and size $y$ gives birth
to a daughter cell of size $x<y$. Integrating over the daughter cell's
size $x$, $D$ and $\tilde{D}$ satisfy $D(a, y, t) =
\int_0^{y}\tilde{D}(a, y, x, t)\dd x$, reflecting cell number
conservation. Finally, to maintain biomass conservation during
division, $\tilde{D}(a, x, y, t) = \tilde{D}(a, x, x-y,
t)$. The prefactor 2 in Eq.~\eqref{BOUNDARY} indicates that a cell of
size $y$ gives birth to one daughter cell of size $y-x$ and
another of size $x$.

The nonlocal boundary condition \eqref{BOUNDARY} for cell
proliferation plays an essential role in depicting how cell division
affects the cell population size and age structure, and presents a
major obstacle in numerical computation as the integration is taken in
the unbounded domain $(x, +\infty)\times(0, +\infty)$. Another numerical
    challenge arises from a possible ``blow-up'' behavior in which
    
\begin{equation}
\lim\limits_{t\to+\infty}\langle x(t)\rangle =
\f{\int_0^{+\infty} \int_0^{+\infty} x n(a,x,t) \dd{a}\dd{x}}
  {\int_0^{+\infty}\int_0^{+\infty} n(a,x,t)\dd{a}\dd{x}} = +\infty.
\label{blowup}
\end{equation}
Whether blowup can occur is of biological interest
\cite{KESSLER2017,Xia2020} and has been predicted within certain cell
proliferation models \eqref{PDEForm} under specific
conditions~\cite{KESSLER2017}.

Existing numerical methods such as the finite volume method in \cite{Xia2020} typically
truncate the unbounded domain into a bounded domain and therefore
cannot accurately capture long time blowup behavior of $\langle
x(t)\rangle$. The need for numerical solutions in the unbounded domain
$\Lambda\times\Lambda$ for Eqs.~\eqref{PDEForm} and \eqref{BOUNDARY}
is thus evident. We apply the scaling technique built in to
Alg.~\ref{algscaling} only in $x$-dimension for tracking the
increasing $\langle x(t)\rangle$, considering the age distribution is
often presumed to be stable since no cell could live too long without
division.  A standard two-dimensional pseudo-spectral method with the
generalized Laguerre functions are used in $(a,x)$-space, coupled with
a third-order TVD Runge-Kutta time discretization in $t$. % \cite{Tang2006Higher}.

\begin{example}
\label{ex:PDE2Dsolve}
\rm 
We solve Eqs.~\eqref{PDEForm} and \eqref{BOUNDARY} with $g(a, x,
t)=t+7$, $\sigma(a, x, t) = 2(t+6)x$, $D(a, x, t) = {x}/{(t+5)}$,
$\tilde{D}(a, y, x, t) = {1}/{(t+5)}$. These parameters leads to the
analytic solution $n(a, x, t) = e^{t}e^{-2a} \exp(-x/(5+t))$, which
produces the mean size $\langle x(t)\rangle = 5+t$. This result shows
that the average size is unbounded as it grows linearly in time and
thus, for general cases, requires proper scaling in $x$-dimension.  We
adopt the same expansion order $N$ in both size $x$- and age $a$-dimensions.  For
the nonlocal boundary condition given in Eq.~\eqref{BOUNDARY}, we also
use $N+1$ Laguerre-Robatto collocation points in each dimension to
perform the numerical integration.

Fig.~\ref{fig7} presents the numerical results with the initial
scaling factors $(\beta_a, \beta_x)=(1, 0.9)$ and a timestep of
$0.002$.  We observe that the frequency-dependent scaling
algorithm~\ref{algscaling} in $x$-dimension shows a faster spectral
convergence with $N$ than that of the unscaled algorithm (see
Fig.~\ref{fig7}(a)).  That is, both the sizer-timer model
\eqref{PDEForm} in unbounded domain and the nonlocal boundary
condition \eqref{BOUNDARY} are well resolved by the Laguerre spectral
approximation with frequency-dependent scaling.  When fixing $N=20$,
the unscaled numerical solution experiences an error growth to
$1.143\text{e-}02$ till $t=10$ for using inappropriate scaling
factors, whereas the error of the scaled solution is less than
$8.662\text{e-}06$ (see Fig.~\ref{fig7}(b)). The frequency indicator
in the $x$-dimension is kept around $10^{-6}$ (red curve with
left-pointing triangles in Fig.~\ref{fig7}(c)) by continuously
shrinking the scaling factor $\beta_x$ from $0.9$ to $0.2766$ for
tracking the blowup (black curve with asterisks in
Fig.~\ref{fig7}(d)).  The average size of the scaled solution behaves
almost exactly like $\langle x(t)\rangle = 5+t$ and the value at
$t=10$ is $15.001$ (see red curve with left-pointing triangles in
Fig.~\ref{fig7}(d)).  Note that the scaling in $a$-dimension will
really not be triggered even when we apply the scaling algorithm for
both $x$- and $a$-dimensions.

\begin{figure}[htb]
        \begin{center}
    \includegraphics[width=4.4in]{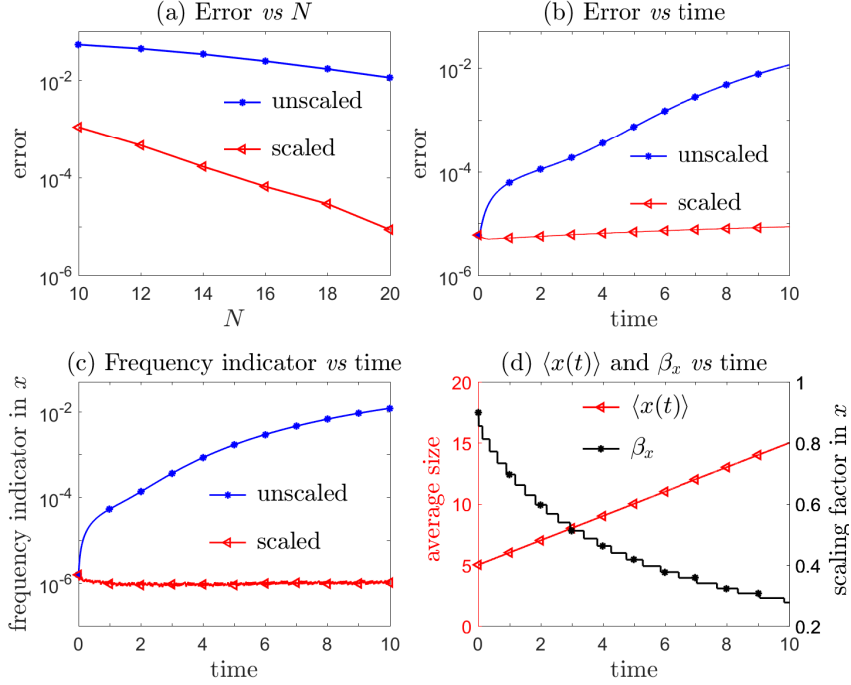}
	\end{center}
\caption{\small Numerical results obtained by the scaling
  algorithm~\ref{algscaling} for the structured cell population
  proliferation model~\eqref{PDEForm} with the nonlocal
  boundary~\eqref{BOUNDARY}: The scaled method gives better results
  than the unscaled one till $t=10$. The latter experiences a growth
  in error because inappropriate scaling factors are used, whereas the
  former gains a faster spectral convergence in the expansion order
  $N$.  We adopt the same $N$ in both size $x$- and age $a$-dimensions
  and set $N=20$ for the last three plots. The frequency-dependent
  scaling is applied only in $x$-dimension for tracking the blowup
  behavior in Eq.~\eqref{blowup}. The frequency indicator in
  $x$-dimension is kept around $10^{-6}$ through constantly shrinking
  the scaling factor $\beta_x$ to capture the blowup.  The average
  size of the scaled solution is in good agreement with that of the
  analytical solution, \textit{i.e.}, $\langle x(t)\rangle=5+t$.}
     \label{fig7}
\end{figure}
\end{example}

\section{Summary and Conclusions}
\label{sec:con}

The key to making spectral approximations in unbounded domains more
efficient is to allocate collocation points in an economical manner
such that crucial regimes of unknown solutions can be resolved
accurately. This is essentially an adaptive numerical method for PDEs
in unbounded domains, for which there are very few studies compared
with its bounded-domain counterpart. Using the standard language of
adaptive methods, the proposed scaling technique based on the
frequency indicator can be regarded as $r$-adaptivity to redistribute
collocation points via adjusting the scaling factor, while the
proposed moving technique based on the exterior-error indicator is
similar to $h$-adaptivity to add collocation points in the interior
subdomain.  Both indicators utilize only the numerical solution and do
not require any a prior knowledge of unknown solutions. The frequency
indicator can be also used in a refinement technique \cite{XIA2020b} which
corresponds to the $p$-adaptivity, useful for time-dependent problems
with oscillations at infinity.

\section*{Acknowledgments}

MX and TC acknowledge support from the National Science Foundation
through grant DMS-1814364 and the Army Research Office through grant
W911NF-18-1-0345.  SS acknowledges the financial support from the
National Natural Science Foundation of China (Nos.~11822102,
11421101), Beijing Academy of Artificial Intelligence (BAAI) and the
computational resource provided by High-performance Computing Platform
of Peking University.

\bibliographystyle{siamplain}
%\bibliography{bibliography}

\end{document}